\documentclass[twocolumn,twoside,10pt]{IEEEtran}

\usepackage{amsfonts,amsmath,graphicx,algorithmic}

\begin{document}

\title{An exact tree projection algorithm for wavelets}
\author{Coralia Cartis and Andrew Thompson*
\thanks{C. Cartis is with the School of Mathematics and the Maxwell Institute, University of Edinburgh, Edinburgh, UK (coralia.cartis@ed.ac.uk).}
\thanks{A. Thompson is with the Department of Mathematics, Duke University, Durham, North Carolina (thompson@math.duke.edu).}}

\maketitle

\newcommand{\argmax}{\operatornamewithlimits{argmax}}
\newcommand{\argmin}{\operatornamewithlimits{argmin}}
\newcommand{\mO}{\mathcal{O}}
\newcommand{\mP}{\mathcal{P}}
\newcommand{\mC}{\mathcal{C}}
\newcommand{\Tk}{\mathcal{T}_k}
\newcommand{\RR}{\mathbb{R}}
\newcommand{\ZZ}{\mathbb{Z}}
\newcommand{\NN}{\mathbb{N}}
\newcommand{\Gm}{\Gamma}
\newcommand{\sh}{\hat{s}}
\newcommand{\GG}{G}
\newcommand{\tGG}{\tilde{G}}
\newcommand{\tF}{\tilde{F}}
\newcommand{\Om}{\Omega}

\newcommand{\comment}[1]{}

\newtheorem{thm}{Theorem}
\newtheorem{lem}[thm]{Lemma}

\begin{abstract}
We propose a dynamic programming algorithm for projection onto wavelet tree structures. In contrast to other recently proposed algorithms which only give approximate tree projections for a given sparsity, our algorithm is guaranteed to calculate the projection exactly. We also prove that our algorithm has $\mO(Nk)$ complexity, where $N$ is the signal dimension and $k$ is the sparsity of the tree approximation.
\end{abstract}

\begin{IEEEkeywords}
Sparse representations, wavelets, dynamic programming, complexity analysis, compressed sensing.\\ 
EDICS categories: DSP-SPARSE, MLSAS-SPARSE.
\end{IEEEkeywords}

\section{Introduction}\label{intro}

The discrete wavelet transform is a much-used tool in digital signal processing, especially in image processing, since it provides sparse representations for piecewise-smooth signals~\cite{mallat}. Wavelets have a multi-scale structure in which a signal is convolved with dilations and translations of some `mother wavelet', which induces a natural dyadic tree structure upon the wavelet coefficients. Since wavelets essentially detect discontinuities, large wavelet coefficients of piecewise smooth signals are often propagated down branches of the tree. This suggests that sparse approximations of such signals have additional structure: they are \textit{tree-sparse}, by which we mean that they are supported on a rooted subtree~\cite{modelbased}, see Section~\ref{prelim}. This paper concerns the task of finding an exact tree projection, namely a tree-sparse vector of a given sparsity which is closest in Euclidean norm to some given signal vector. 

Tree projections, or approximations thereof, predate the arrival of wavelets. An early example of their computation is found in the pruning of decision tree structures known as classification and regression trees (CART)~\cite{breiman,bohanec_bratko,cartbob}. With the emergence of wavelets, tree approximations also began to be used in the context of wavelet-based compression~\cite{shapiro,cohen} and orthogonal best-basis methods~\cite{coifman,cartbob,lagrangian2}. More recently, algorithms for tree-based Compressed Sensing have been proposed which require the computation of tree projections~\cite{union_subspaces,modelbased}. Furthermore, recovery guarantees have been obtained for these algorithms which depend crucially on the assumption that an exact tree projection is computed in every iteration of the algorithm~\cite{union_subspaces,modelbased}.

Several algorithms have been proposed for wavelet tree projection, and the summary  by Baraniuk in~\cite{treeapprox} describes three possible approaches: \textit{Greedy Tree Approximation} (GTA), the \textit{Condensing Sort and Select Algorithm} (CSSA), and the \textit{Complexity-Penalized Sum of Squares} (CPSS). GTA~\cite{cohen} proceeds by growing the tree from the root node in a `greedy' fashion by making a series of locally optimal choices. However, such an approach is only guaranteed to find an exact tree projection in the idealized case in which the wavelet coefficients are monotonically nonincreasing along the branches~\cite{treeapprox}. The CSSA was originally proposed by Baraniuk and Jones~\cite{CSSA} in the context of optimal kernel design, and solves a linear programming (LP) relaxation of the exact tree projection problem. The CPSS approach has its roots in the \textit{minimal cost-complexity pruning} algorithm for CART~\cite{breiman}, and was later developed by Donoho in the context of wavelet approximations~\cite{cartbob,lagrangian2}. This algorithm solves a Lagrangian relaxation of the exact tree projection problem, in which the sparsity constraint is removed and penalized in the objective function instead. We argue in Section~\ref{prelim} that, since both the CSSA and CPSS solve relaxations, neither algorithm is guaranteed to find an exact tree projection for a given sparsity. Since the theoretical results mentioned above require an exact tree projection, a crucial question is whether it is possible to guarantee the computation of an exact tree projection in polynomial time. 

In this paper, we propose an algorithm for exact tree projection in the context of wavelets. We make use of a Dynamic Programming (DP) approach which has been used before in other contexts. An early reference is Breiman et al.~\cite{breiman}, where it is proved that such an algorithm exists in the context of CART decision trees. The algorithm, referred to as \textit{optimal tree pruning}, was then formally stated by Bohanec and Bratko~\cite{bohanec_bratko}, who also proved that it has complexity $\mO(N^2)$, where $N$ is the signal dimension. The same exact tree projection approach was also used in the context of orthogonal best-basis selection in~\cite{coifman}. Our algorithm adopts the same basic approach as in~\cite{bohanec_bratko}, but this time in the context of wavelet approximation; further distinctions to~\cite{bohanec_bratko} are given in Section~\ref{our_algorithm}. We also prove that our algorithm has low-order polynomial complexity, namely $\mO(Nk)$, where $k$ is the sparsity of the tree approximation, which represents an improvement over the result in~\cite{bohanec_bratko}.
 
Some closely-related recent work~\cite{pairwise} proves that there exists a polynomial-time DP algorithm for projection onto any tree structure which has a \textit{loopless pairwise overlapping group} property. We should point out, however, that wavelet transforms may not exhibit this property.

\section{Further motivation and preliminaries}\label{prelim}

\begin{figure}[h]
\centering
\includegraphics[width=2.5in,height=1.8in]{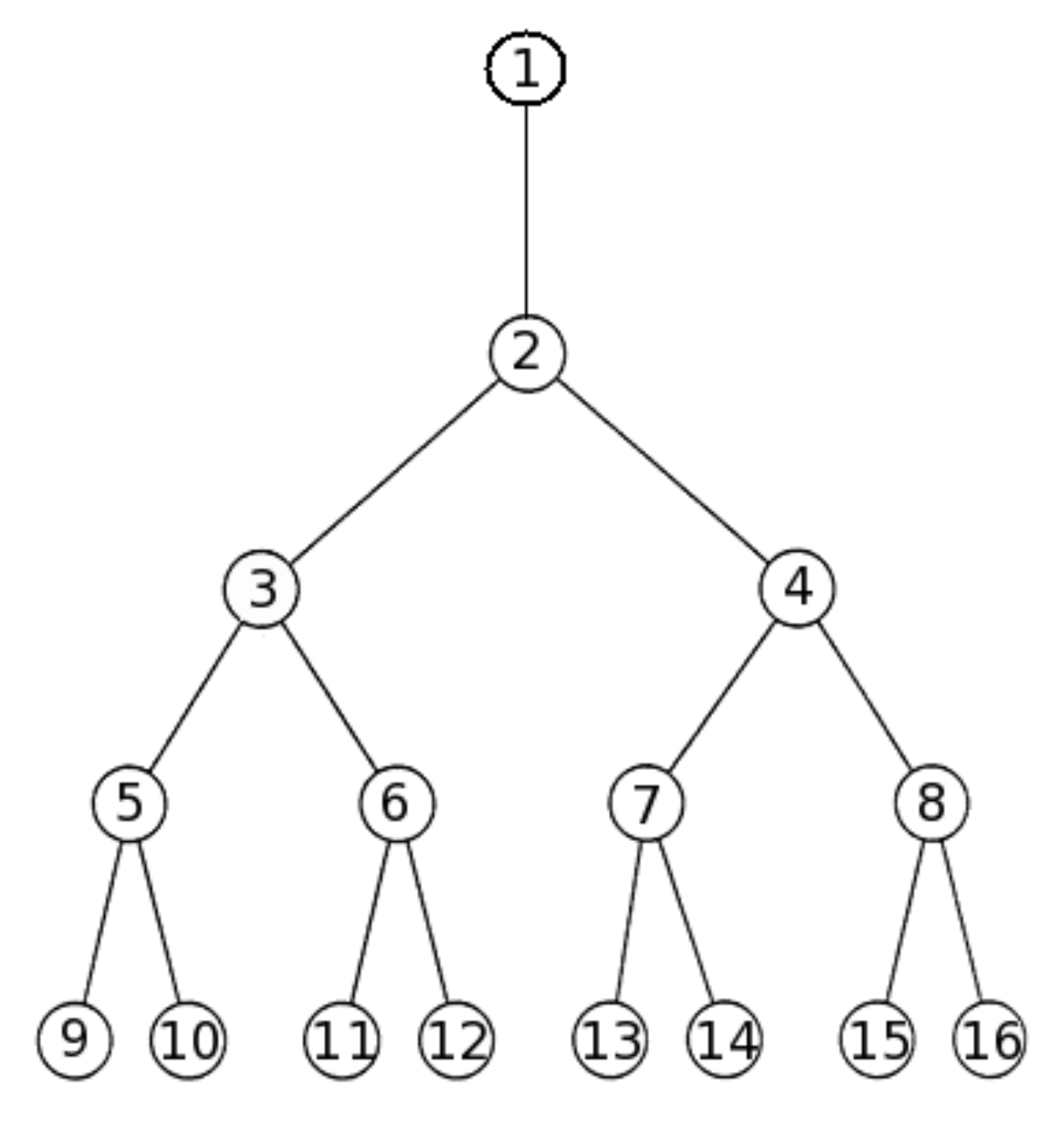}
\caption{An illustration of the canonical dyadic tree structure considered here, for the case of $d=2$ (binary tree) and $J=4$.}
\label{tree_number}
\end{figure}

We frame the discussion in terms of standard $D$-dimensional Cartesian-product dyadic wavelets, which have a particular canonical $2^D$-ary tree structure. We assume a tree structure in which each coefficient has a maximum of $d$ children, where $d$ is the tree order. The Cartesian-product structure implies that $d=2^D$, where the transform dimension $D$ is some fixed positive integer, and so $d$ is a fixed integer with $d\geq 2$. Also, let $N=d^J$ for some $J\geq 2$, where $J$ is the number of levels in the tree structure. We assign a numbering $\{i:1\le i\le N\}$ to the nodes of the tree, starting with the root node and proceeding from coarse to fine levels. Let the root node $i=1$ be at level $0$, and have $d-1$ children at level $1$, namely $\{2,\ldots,d\}$. Let each of the remaining nodes except those in level $J$ (the finest level), that is $\{i:2\le i\le d^{J-1}=N/d\}$, each have $d$ children, namely $\{d(i-1)+1,\ldots,di\}$ respectively. The nodes in the finest level $J$, that is $\{i:i>d^{J-1}=N/d\}$, are assumed to have no children; see Figure~\ref{tree_number}.

Following~\cite{modelbased}, we define a \textit{rooted tree of cardinality} $k$ to be some subtree consisting of $k$ nodes of the full tree just described, subject to the restrictions that it must include the root node, and that if any node other than the root node is included then its parent node must also be included. We write $\Tk$ for the set of supports of all rooted trees of cardinality $k$. Given a vector $y\in\RR^N$, we are interested in finding the Euclidean projection $\mP_k(y)$ of a vector $y\in\RR^N$ onto the set of vectors with support in $\Tk$, namely
\begin{equation}\label{treeproj_def}
\mP_k(y)=\argmin_{\{z\in\RR^N:\mathrm{supp}(z)\in\Tk\}}\|z-y\|_2,
\end{equation}
where $\mathrm{supp}(z)$ denotes the support of $z$. Provided the above numbering scheme is used to distinguish between supports with precisely the same approximation error, $\mP_k(y)$ is well-defined. Intuitively, the following lemma establishes that a tree projection preserves the values of selected coefficients. 
\begin{lem}\label{Pk_proj}
Let $\mP_k(\cdot)$ be defined as in (\ref{treeproj_def}). Then
\begin{equation}\label{treeproj_def2}
\left\{\mP_k(y)\right\}_i:=\left\{\begin{array}{ll}
y_i&i\in\Gm\\
0&i\notin\Gm
\end{array}\right.\;\;\;\;
\mbox{where}\;\;\;\;\Gm:=\argmax_{\Omega\in\Tk}\left\|y_{\Omega}\right\|.
\end{equation}
\end{lem}

\proof{Let $\mbox{supp}(z)=\Gm\in\Tk$. Then
$$\|z-y\|^2=\left\|(z-y)_{\Gm}\right\|^2+\left\|y_{\Gm^C}\right\|^2,$$
and $\left\|(z-y)_{\Gm}\right\|^2$ is minimized by setting $z_i=y_i$ for all $i\in\Gm$. Meanwhile, $\left\|y_{\Gm^C}\right\|$ is minimized by choosing $\Gm$ to be the set $\Omega\in\Tk$ which maximizes $\left\|y_{\Omega}\right\|$, and the result now follows.\hfill$\Box$}

Lemma~\ref{Pk_proj} implies that tree projection is a decision problem: since coefficient values are preserved, the projection is obtained by identifying the set $\Gm\in\Tk$ on which $y$ has maximum energy. It follows that 
$\mP_k(y)$ may also be found by solving the following Integer Program (IP):
\begin{equation}\label{treeIP}\displaystyle\max_{\tau\in\ZZ^N}\sum_{i=1}^N y_i^2 \tau_i\;\;\;\mbox{subject to}\;\;\;\displaystyle\begin{array}{l}\tau\geq 0\\
\{\tau_i\}\;\mbox{tree-nonincreasing}\\
\sum\tau_i=k\\
\tau_1=1.\end{array}
\end{equation}
Here tree-nonincreasing means that the coefficients do not increase along the branches, a condition which may be translated into a series of linear constraints. Note that the assumption of integrality, together with the constraints, in fact forces $\tau_i\in\{0,1\}$ for each $i$, so that $\tau$ acts as a decision variable. Denoting the optimal solution of (\ref{treeIP}) by $\tau^{\ast}$, the coefficients of the best tree approximation are then $y_i \tau^{\ast}_i$. 

We can now be more precise about our claim in Section~\ref{intro} that both the CSSA~\cite{CSSA} and the CPSS~\cite{cartbob} solve relaxations of the tree projection problem, and hence are only guaranteed to give approximate tree projections. 

The CSSA solves an LP relaxation of (\ref{treeIP}) in which one dispenses with the assumption of integrality while retaining all other constraints. The authors of~\cite{CSSA} observe that the two solutions are sometimes equal, but not always: a value of $\tau_i=1$ may be assigned to all but a few coefficients, which are each assigned a value $\bar{\tau}$ for some $0<\bar{\tau}<1$. In this case, the result is a tree of size strictly greater than $k$. Moreover, there is no straightforward method for `adjusting' the solution \textit{a posteriori} to obtain the optimal $k$-sparse tree projection. 

Meanwhile, the CPSS~\cite{cartbob} solves a Lagrangian relaxation of (\ref{treeIP}) in which the equality constraint $\sum\tau_i=k$ is removed and penalized in the objective instead. This reformulation can be solved using fine-to-coarse dynamic programming on the tree~\cite{cartbob,lagrangian2}. Because (\ref{treeIP}) is an IP, its Lagrangian relaxation is also not guaranteed to share the same optimal solution as (\ref{treeIP})~\cite{integer}. A further drawback of the CPSS is the non-obvious relationship between the Lagrange multiplier $\lambda$ and the required sparsity, which means that the problem must in practice be solved for multiple values of $\lambda$ in order to find a tree approximation for a given sparsity.

\section{The ETP algorithm}\label{our_algorithm}

Our algorithm (Algorithm 1) falls into the broad category of dynamic programming (DP) algorithms which optimize on directed graphs by utilizing a principle of optimality, namely that optimal solutions at a given node may be determined entirely from optimal solutions at `preceding' nodes~\cite{bellman}. Our algorithm makes two passes through the tree: firstly, a fine-to-coarse pass finds optimal subtrees at each node for all $\tilde{k}\le k$. Secondly, a coarse-to-fine pass tracks back to identify the optimal solution.

The algorithm starts at the finest level and moves up the tree, finding optimal subtrees rooted at each node, each decision requiring only the information already obtained within that particular subtree. Following~\cite{bohanec_bratko}, optimal subtrees for a given node are obtained by successively incorporating each of its children, and repeatedly optimizing over all known solutions. To maximize efficiency, optimal subtrees are in fact only found for all cardinalities which could possibly contribute to a rooted tree of cardinality $k$, which imposes two restrictions. Firstly, by summing the appropriate geometric series, the maximum cardinality of a subtree rooted at a node in level $j$ is $\frac{d^{J+1-j}-1}{d-1}$. Secondly, any subtree rooted at a node in level $j$ of cardinality greater than $k-j$ would necessarily contribute to a rooted tree of cardinality greater than $k$. Defining
\begin{equation}\label{lj_def}
l(j):=\min\left(\frac{d^{J+1-j}-1}{d-1},k-j\right)\;\;1\le j\le J,
\end{equation}
we therefore need only find optimal subtrees for cardinalities up to $l(j)$ for nodes at level $j$. Further restrictions are imposed due to the fact that children are incorporated sequentially.

$F(i,l)$ denotes the total energy of the optimal subtree rooted at node $i$, of size $l$, and $\GG(i,l)=[G_r(i,l):r=1,\ldots,d]\in\RR^d$ is a vector of nonnegative integers, giving the number of nodes contributing to this optimal solution from the subtree rooted at each of the children of node $i$. We then proceed recursively and eventually determine $F(1,k)$, the energy of the optimal tree of size $k$.\footnote{We also need to define two temporary variables $\tilde{F}$ and $\tilde{G}$ which represent updates to $F$ and $G$ during the incorporation of a child.} Finally, we use the precedence information to trace the optimal solution back along the branches. The algorithm actually does more than is asked for: by the time the root node is reached at the end of the first pass, enough information has been obtained to determine optimal trees for all $\tilde{k}\le k$. 

ETP is closely related, but not identical, to the optimal tree pruning algorithm in~\cite{bohanec_bratko}, where the decision tree context motivates the solution of a slightly different problem from the present tree projection problem. In their setup, if a node is pruned, then so are all other nodes sharing the same parent in the tree structure. Also, an optimal pruning is obtained for a given number of \textit{leaves} (ends of branches of the tree), as opposed to the cardinality of the tree approximant.  There is also a difference in how information is stored: the proposal in~\cite{bohanec_bratko} is to store the full optimal subtrees themselves at each node, essentially incorporating the backtracking phase into the first pass of the algorithm.\\
\\
\textbf{Algorithm 1. Exact tree projection (ETP).}\\
\textbf{Inputs:} $\displaystyle\left\{\begin{array}{ll}
d\in\NN&(d\geq 2)\\
y\in\RR^N&(N=d^J;\;J\geq 2)\\
k\in\NN&(2\le k\le N)
\end{array}\right.$\\
\textbf{Find all optimal subtrees:}
\begin{algorithmic}
\FOR {$i=(d^{J-1}+1):d^J$}
\STATE $F(i,0)=0$ and $F(i,1) = y^2_i$
\ENDFOR
\FOR {$j=(J-1):-1:1$}
\FOR {$i=(d^{j-1}+1):d^j$}
\STATE $F(i,0)=0$ and $F(i,1)=y^2_i$
\STATE $\GG(i,0)=0$ and $\GG(i,1)=0$
\FOR {$r=1:d$}
\FOR {$l=2:\min[l(j),rl(j+1)+1]$}
\STATE $s_-=\max\{1,l-[(r-1)l(j+1)+1]\}$
\STATE $s_+=\min[l-1,l(j+1)]$
\STATE $\sh=\displaystyle\argmax_{s_-\le s\le s_+}\left\{F[d(i-1)+r,s]+F(i,l-s)\right\}$
\STATE $\tF(i,l)=F[d(i-1)+r,\sh]+F(i,l-\sh)$
\STATE $\tGG(i,l)=\GG(i,l-\sh)$, $\tGG_r(i,l)=\sh$
\ENDFOR
\FOR {$l=2:\min[l(j),rl(j+1)+1]$}
\STATE $F(i,l)=\tF(i,l)$ and $\GG(i,l)=\tGG(i,l)$
\ENDFOR
\ENDFOR
\ENDFOR
\ENDFOR
\FOR {$r=2:d$}
\FOR {$l=2:\min[k,(r-1)l(1)+1]$}
\STATE $s_-=\max\{1,l-[(r-2)l(1)+1]\}$
\STATE $s_+=\min[l-1,l(1)]$
\STATE $\sh=\displaystyle\argmax_{s_-\le s\le s_+}\left[F(r,s)+F(1,l-s)\right]$
\STATE $\tF(1,l)=F(r,\sh)+F(1,l-\sh)$
\STATE $\tGG(1,l)=G(1,l-\sh)$, $\tGG_r(1,l)=\sh$
\ENDFOR
\FOR {$l=2:\min[k,(r-1)l(1)+1]$}
\STATE $F(1,l)=\tF(1,l)$ and $\GG(1,l)=\tGG(1,l)$
\ENDFOR
\ENDFOR
\end{algorithmic}
\textbf{Backtrack to identify solution:}
\begin{algorithmic}
\STATE $\tau=0$, $\tau_1=1$, $\Gm_1=k$
\FOR {$j=0:(J-1)$}
\FOR {$i=\max[2,(d^{j-1}+1)]:d^j$}
\IF {$\tau_i=1$}
\FOR {$r=\max(1,2-j):d$}
\IF {$G_r(i,\Gm_i)>0$}
\STATE $\tau_{d(i-1)+r}=1$
\STATE $\Gm_{d(i-1)+r}=G_r(i,\Gm_i)$
\ENDIF
\ENDFOR
\ENDIF
\ENDFOR
\ENDFOR\\
\end{algorithmic}
\textbf{Calculate solution:}
\begin{algorithmic}
\FOR {$i=1:N$}
\STATE $\hat{y_i}=y_i\tau_i$
\ENDFOR
\end{algorithmic}
\textbf{Output:} $\hat{y}\in\RR^N$

\section{Complexity analysis}

To establish the complexity of ETP, we note that the dominating operations are additions and comparisons, and therefore it suffices to bound the total number of these operations. 

\begin{thm}\label{complexitybound}
Given $y\in\RR^N$, the ETP algorithm requires at most $(3d^2 Nk+N)$ additions/comparisons to calculate $\mP_k(y)$, defined in (\ref{treeproj_def}).
\end{thm}

\proof{Consider the first pass of ETP. In each level $j$ with $1\le j\le J-1$, there are $(d-1)\cdot d^{j-1}$ nodes. For a given node $i$ at level $j$, each of its $d$ children are incorporated successively. Each time a child is incorporated, we calculate optimal subtrees for each $l$ such that $2\le l\le l_{max}$. Each of these calculations requires at most $(l-1)$ additions, each one accompanied by a comparison, and we have $l_{max}\le l(j)$, where $l(j)$ is defined in (\ref{lj_def}). The number of additions/comparisons needed to incorporate a child is therefore bounded by
\begin{equation}\label{child}
2\sum_{l=2}^{l(j)}(l-1)\le l(j)[l(j)-1]\le[l(j)]^2.
\end{equation}
Combining all these observations, and writing $\mC$ for the total number of additions/comparisons in the first pass of the algorithm for levels $1\le j\le J-1$ (i.e. excluding the root node), we have
\begin{equation}\label{firstbound}
\mC\le \displaystyle\sum_{j=1}^{J-1}\left\{(d-1)d^{j-1}\cdot d\cdot[l(j)]^2\right\},
\end{equation}
To bound $l(j)$, observe from (\ref{lj_def}) that
\begin{equation}\label{lmaxbound}
l(j)\le\min(d^{J+1-j},k).
\end{equation}
To determine which of $d^{J+1-j}$ or $k$ gives a tighter bound, let $1\le p\le J-1$ be the unique positive integer such that 
\begin{equation}\label{pbounds}
d^{p-1}\le k\le d^p.
\end{equation} 
First let us assume $p>2$. We have
$$k\le d^{J+1-j}\Longleftrightarrow d^p\le d^{J+1-j}\Longleftrightarrow j\le J+1-p,$$
which we may combine with (\ref{firstbound}) and (\ref{lmaxbound}) to deduce
\begin{equation}\label{middlebound}
\mC\le\displaystyle\sum_{j=1}^{J+1-p}(d-1)d^j\cdot k^2+\displaystyle\sum_{j=J+2-p}^{J-1}(d-1)d^j\cdot d^{2(J+1-j)}
\end{equation}
$$=(d-1)\left\{dk^2\displaystyle\sum_{j=1}^{J+1-p}d^{j-1}+\displaystyle\sum_{j=J+2-p}^{J-1}d^{2J+2-j}\right\}$$
\begin{equation}\label{secondbound}
=(d-1)\left\{dk^2\left(\frac{d^{J+1-p}-1}{d-1}\right)+d^{J+p}\left[\frac{1-\left(\frac{1}{d}\right)^{p-2}}{1-\frac{1}{d}}\right]\right\}.
\end{equation}
Since 
$$(d-1)\left[\displaystyle\frac{1-\left(\frac{1}{d}\right)^{p-2}}{1-\left(\frac{1}{d}\right)}\right]\le\left(\frac{d-1}{1-\frac{1}{d}}\right)=d,$$
and since $d^{J+1-p}-1\le d^{J+1-p}$, we can further deduce from (\ref{secondbound}) that
$$\mC\le k^2\cdot d^{J+2-p}+d^{J+1+p},$$
to which we can make the substitution $N=d^J$ and apply the bounds (\ref{pbounds}), concluding that
\begin{equation}\label{bound1}
\mC\le k^2\cdot\displaystyle\frac{d^2 N}{k}+d^2 Nk=2d^2 Nk.
\end{equation}
On the other hand, if $p\le 2$, then $k\le d^2\le d^{J+1-j}$ for all $j\le J-1$, and so the second summation in (\ref{middlebound}) is empty. We may then follow the same argument to bound the first summation, obtaining $\mC\le 2d^2 Nk$ in this case also. Meanwhile, the root node has $(d-1)$ children, which may be combined with (\ref{child}) and (\ref{lmaxbound}) to deduce that the number of additions/comparisons for the root node is bounded by 
\begin{equation}\label{bound2}
(d-1)k^2\le d^2 Nk.
\end{equation}
Finally, note that the second pass involves no additions and at most $N$ comparisons, which combines with (\ref{bound1}) and (\ref{bound2}) to yield the result.\hfill$\Box$}

It follows from Theorem~\ref{complexitybound} that ETP has complexity $\mO(Nk)$.\footnote{Note that, except for the constant involved, the result is independent of $d$.} We may compare the complexity of ETP with that of the other approximate tree projection methods discussed in Section~\ref{prelim}: GTA~\cite{cohen} and the CSSA~\cite{CSSA} are both $\mathcal{O}(N\log N)$ while the CPSS~\cite{cartbob} is $\mathcal{O}(N)$. Provided $\log N\ll k$, we can see that the price we pay for guaranteeing an exact tree projection is an increased order of complexity.

The algorithm is easy to implement and we have experimented with random tests successfully; its cost is not disappointing compared to approximate methods such as the CSSA and the CPSS.\footnote{Matlab code is available to download from \texttt{www.math.duke.edu/\~{}thompson}.}

\section{Conclusion}\label{conclusion}

We have presented an algorithm for exact tree projection in the context of wavelets and have proved that it has $\mO(Nk)$ worst-case complexity. The DP algorithm described here could also be applied to any tree-based signal representations. Approximate algorithms may still offer an advantage in terms of computational efficiency, and therefore the practitioner's choice of tree projection algorithm will depend upon whether guaranteed precision or algorithm running time is the overriding priority.

\bibliographystyle{IEEEtran}
\bibliography{treeprojbib}

\comment{
\FOR {$r=2:d$}
\FOR {$l=2:min[k,(r-1).l(1)+1]$}
\STATE $s_{min}=max\{1,l-min[k,(r-2).l(1)+1]\}$
\STATE $s_{max}=min[l-1,l(1)]$
\STATE $\sh=\displaystyle\argmax_{s_{min}\le s\le s_{max}}\left[F(r,s)+F(1,l-s)\right]$
\STATE $F(1,l)=F(r,\sh)+F(1,l-\sh)$
\STATE $\GG(1,l)=G(1,l-\sh)$
\STATE $G_r(1,l)=\sh$
\ENDFOR
\FOR {$l=2:min[k,(r-1).l(1)+1]$}
\STATE $F(1,l)=\tF(1,l)$ and $\GG(1,l)=\tGG(1,l)$
\ENDFOR
\ENDFOR
}

\end{document}